\documentclass[a4paper,12pt]{amsart}   

\usepackage[T1]{fontenc}
\usepackage{amsmath,amssymb,amsxtra,amscd}
\usepackage{amsthm,NewStableClassification}

\newtheorem{theorem}{Theorem}[section]
\newtheorem{proposition}[theorem]{Proposition}
\newtheorem{lemma}[theorem]{Lemma}
\newtheorem{definition}[theorem]{Definition}
\newtheorem{example}[theorem]{Example}

\newtheorem*{theorem*}{Theorem}
\newtheorem*{question*}{Question}

\theoremstyle{definition}

\begin{document}

\title{Alternative stable homotopy classification of $\pComp{BG}$}
\author{K\'ari Ragnarsson}
\date{September 12, 2005}
\thanks{The author is supported by EPSRC grant GR/S94667/01}
 \address{Department of Mathematical Sciences, University of Aberdeen, Aberdeen AB24 3UE, United Kingdom.}
\email{kari@maths.abdn.ac.uk}
\maketitle

\begin{abstract}
We give an alternative to the stable classification of
$p$-completed homotopy types of classifying spaces of finite
groups offered by Martino-Priddy in \cite{MP3}. For a finite group
$G$ with Sylow subgroup $S$, we regard the stable $p$-completed
classifying space $\Stable{\pComp{BG}}$ as an object under
$\Stable{BS}$ via the canonical inclusion map. Thus we get a
classification in terms of induced fusion systems. Applying
Oliver's solution \cite{Ol1,Ol2} to the Martino-Priddy conjecture
\cite{MP}, we obtain the surprising result that the unstable
homotopy type of $\pComp{BG}$ is determined by the map \mbox{$
\Stable{BS} \to \Stable{\pComp{BG}}$}, but not by the homotopy
type of $\Stable{\pComp{BG}}$.
\end{abstract}

\section*{Introduction}
Using Carlsson's solution of the Segal Conjecture \cite{Car}, and
Nishida's theory of dominant summands \cite{Ni}, Martino-Priddy
proved a classification theorem for stable homotopy types of
$p$-completed classifying space of finite groups in \cite{MP3}.
For a set $X$, let $\Fp X$ denote the $\Fp$-vector space with
basis $X$ and for groups $Q$ and $G$, let $\Rep{Q}{G} =
\Hom{Q}{G}{}/G$ with $G$ acting by conjugation, and let
$\InjRep{Q}{G} \subset \Rep{Q}{G}$ be the set of conjugacy classes
of  injective homomorphisms.
\begin{theorem*}[Stable classification \cite{MP3}]\label{thm:MPStableClassification}
For two finite groups $G$ and $G'$, the following are equivalent:
\begin{itemize}
  \item[(1)] $\pComp{BG}$ and $\pComp{BG'}$ are stably homotopy
  equivalent.
  \item[(2)] For every $p$-group $Q$,
  $$\Fp \Rep{Q}{G} \cong \Fp \Rep{Q}{G'} $$
  as $\Out{Q}{}$-modules.
  \item[(3)] For every $p$-group $Q$,
  $$\Fp \InjRep{Q}{G} \cong \Fp \InjRep{Q}{G'}$$
  as $\Out{Q}{}$-modules.
\end{itemize}
\end{theorem*}

In this paper we take a different point of view. By regarding the
stable $p$-completed classifying space
\mbox{$\ClSpectrum{\pComp{G}} = \Stable{\pComp{BG}}$} of a finite
group $G$ as an object under the stable classifying space
$\ClSpectrum{S}$ of its Sylow subgroup $S$ and using the Segal
conjecture, we are able to reconstruct the fusion system
$\F_S(G)$. This is the category whose objects are the subgroups of
$S$, and whose morphisms are the morphisms induced by conjugation
in $G$. Thus we are able to classify the stable homotopy types
\mbox{$\Stable{\pComp{BG}}$} as objects under the stable
classifying spaces of their Sylow subgroups by their induced
fusion systems. By Oliver's solution of the Martino-Priddy
conjecture \cite{Ol1,Ol2}, the unstable homotopy types of
$p$-completed classifying spaces of finite groups are also
classified by their fusion systems, so we get the following
theorem which is proved in Section
\ref{sec:NewStableClassification}, where the notion of an
isomorphism of fusion systems is also explained.
\begin{theorem*} [Alternative stable classification]
\label{thm:NewStableClassification} For two finite groups  $G$ and
$G'$ with Sylow subgroups $S$ and $S'$, respectively, the
following are equivalent:
\begin{itemize}
  \item[(i)] There is an isomorphism \mbox{$\gamma\negmedspace: S \to
  S'$} and a homotopy equivalence \mbox{$h\negmedspace: \pComp{\ClSpectrum{G}} \to
  \pComp{\ClSpectrum{G'}}$} such that the following diagram commutes up
  to  homotopy:
  \[
\begin{CD}
{\ClSpectrum{S}} @ > {\ClSpectrum{\iota_{S}}} >> \pComp{\ClSpectrum{G}} \\
@ V {\ClSpectrum{\gamma}} VV @ VV h V \\
{\ClSpectrum{S'}} @> {\ClSpectrum{\iota_{S'}}} >> \pComp{\ClSpectrum{G'}}. \\
\end{CD}
\]
  \item[(ii)] There is an isomorphism of fusion systems
  $$(S,\F_{S}(G)) \longrightarrow (S',\F_{S'}(G')).$$
  \item[(iii)] There is a homotopy equivalence of $p$-completed
  classifying spaces
  $$h \negmedspace : \pComp{BG} \longrightarrow \pComp{BG'}.$$
\end{itemize}
\end{theorem*}
The presence of the third condition in this theorem is rather
surprising, since Martino-Priddy have in \cite[Example 5.2]{MP3}
constructed two groups whose $p$-completed classifying spaces are
homotopy equivalent stably, but not unstably. The added datum of
the map \mbox{$\ClSpectrum{S} \to \ClSpectrum{G}$} therefore
really gives new algebraic information. The cost of this added
information is to only allow maps \mbox{$\ClSpectrum{S} \to
\ClSpectrum{S'}$} induced by group homomorphisms. This is an
interesting point, which will be taken up in Section
\ref{sec:Discussion}.

This paper is split into four sections. In Section \ref{sec:Segal}
we recall some background material on the Segal conjecture and in
Section \ref{sec:HoMonos} we do the same for Sylow subgroups of
spaces. In Section \ref{sec:NewStableClassification} we prove the
alternative stable classification theorem and in Section
\ref{sec:Discussion} we compare the two stable classification
theorems to each other as well as to the Martino-Priddy
conjecture.

Throughout this paper, $p$ is a fixed prime. We denote the
inclusion of a group $H$ into a supergroup by $\iota_H$,
specifying the target group if it is not clear from the context.
For a space or a spectrum $X$, let \pComp{X} denote the
Bousfield-Kan $p$-completion of $X$ \cite{Monster}. For spaces $X$
and $Y$, let $\StableMaps{X}{Y}$ denote the group of homotopy
classes of stable maps \mbox{$\Stable{X} \to \Stable{Y}$}. Unless
otherwise specified, spaces and homotopies will be assumed to be
unpointed. We use the shorthand notation $\ClSpectrum(-)$ for the
functor $\Stable{B(-)}$. Recall that \mbox{$\Stable{(\pComp{BG})}
\simeq \pComp{(\Stable{BG})}$} for a finite group $G$ and we write
$\pComp{\ClSpectrum{G}}$ without danger of confusion. Finally,
recall that $p$-completion coincides with $p$-localization for
classifying spaces of finite groups, so we could equally well
state all results for $p$-localizations.

\section{Burnside modules and the Segal conjecture} \label{sec:Segal}
For finite groups $G$ and $G'$, we use the term $(G,G')$-biset to
denote a set with a right $G$-action and a free left $G'$-action
such that the two actions commute. Let $\Morita(G,G')$ denote the
set of isomorphism classes of finite $(G,G')$-bisets. The
operation of taking disjoint unions provides an abelian monoid
structure on $\Morita(G,G')$. We refer to the Grothendieck group
completion of $\Morita(G,G')$ as the \emph{Burnside module of $G$
and $G'$} and denote it by $A(G,G')$. Being an abelian group, we
can regard $A(G,G')$ as a $\Z$-module, and as such its structure
is easy to describe.
\begin{definition}
Let $G$ and $G'$ be finite groups. A \emph{$(G,G')$-pair} is a
pair $(H,\varphi)$ consisting of a subgroup \mbox{$H \leq G$} and
a homomorphism \mbox{$\varphi \negmedspace : H \to G'$}. We say
two $(G,G')$-pairs $(H,\varphi)$ and $(H',\varphi')$ are
\emph{conjugate} if there exist elements \mbox{$g \in G$} and
\mbox{$h \in G'$} such that \mbox{$c_g(H) = H'$} and the following
diagram commutes:
\[
\begin{CD}
{H} @ > {\varphi} >> {G'} \\
@ V \cong V c_g V @ VV c_h V \\
{H'} @> {\varphi'} >> {G'.} \\
\end{CD}
\]
\end{definition}
Conjugacy is an equivalence relation on $(G,G')$-pairs and we
denote the conjugacy class of a $(G,G')$-pair $(H,\varphi)$ by
$[H,\varphi]$. From a $(G,G')$-pair $(H,\varphi)$ we construct a
$(G,G')$-biset
$$G' \times_{(H,\psi)} G = (G'\times G)/\sim,$$
with the obvious right $G$-action and left $G'$-action, where the
equivalence relation $\sim$ is given by
$$(x,gy)\sim (x\psi(g),y)$$
for $x \in G'~,y \in G,~ g \in H$. One can check that this
construction gives a bijection from equivalence classes of
$(G,G')$-pairs to isomorphism classes of indecomposable
$(G,G')$-bisets, and by a slight abuse of notation we will denote
the isomorphism class of the biset \mbox{$G' \times_{(H,\psi)} G$}
also by $[H,\varphi]$. The Burnside module $A(G,G')$ can now be
described as a free $\Z$-module with one basis element for each
conjugacy class of $(G,G')$-pairs.

For a space $X$, let $X_+$ denote the pointed space obtained by
adding a disjoint basepoint to $X$, and let $\PtdStable{X}$ denote
the suspension spectrum of $X_+$. There is a homomorphism
$$\alpha \negmedspace : A(G,G') \longrightarrow \PtdStableMaps{BG}{BG'},$$
defined on basis elements by sending $[H,\varphi]$ to the map
$$\ClSpectrum{\varphi} \circ tr_H \negmedspace : \PtdStable{BG} \xrightarrow{tr_H} \PtdStable{BH} \xrightarrow{\PtdStable{B\varphi}} \PtdStable{BG'},$$
where $tr_H$ is the transfer of the inclusion \mbox{$H
\hookrightarrow G$}. The Segal conjecture, which was proved by
Carlsson in \cite{Car}, states that this homomorphism is a
completion with respect to a certain ideal when $G'$ is the
trivial group. Lewis-May-McClure showed in \cite{LMM} that
consequently the same holds for any finite group $G'$.

When the source group $G$ is a $p$-group, this completion can be
described in a particularly simple way; at least after getting rid
of basepoints. Let $\widetilde{A}(G,G')$ denote the module
obtained from $A(G,G')$ by quotienting out all basis elements of
the form $[P,\varphi]$, where $\varphi$ is the trivial
homomorphism. Recalling that \mbox{$\PtdStable{BG} \simeq
\Stable{BG}\vee \SphereSpectrum$}, one can check that $\alpha$
induces a map
$$\widetilde{A}(G,G') \longrightarrow \StableMaps{BG_+}{BG'_+}/\StableMaps{BG_+}{S^0} \cong \StableMaps{BG}{BG'}.$$
\begin{theorem}[Segal Conjecture \cite{Car,LMM}]\label{thm:Segal}
Let $S$ be a finite $p$-group and $G$ be any finite group. Then
the homomorphism $\alpha$ described above induces an isomorphism
$$\pComp{\widetilde{\alpha}} \negmedspace : \pComp{\widetilde{A}(S,G)} \stackrel{\cong}{\longrightarrow} \StableMaps{BS}{BG},$$
where \mbox{$\pComp{(-)} =(-)\otimes\Zp$} is $p$-adic completion.
\end{theorem}

\section{Homotopy monomorphisms and subgroups} \label{sec:HoMonos}
In this section we describe a homotopy theoretic analogue of group
monomorphisms and subgroup inclusions. The discussion is very
goal-oriented and we develop just enough tools to work with Sylow
subgroups of spaces, as defined below. The reader should be aware
that there are more than one version of homotopy monomorphisms to
be found in the literature and we have here selected one that is
suitable for our purpose.

\begin{definition} A map \mbox{$f: X \to Y$} of spaces is a
\emph{homotopy  monomorphism at $p$} if the induced map in
$\Fp$-cohomology makes $\Coh{X;\Fp}$ a finitely generated
$\Coh{Y;\Fp}$-module.
\end{definition}

As the prime $p$ is fixed we will just say "homotopy
monomorphisms" without danger of confusion. A motivation for the
above definition is the following lemma, a proof of which can be
found in \cite{Ev}.
\begin{lemma} Let $\varphi \negmedspace : P \to Q$ be a
homomorphism of finite $p$-groups. Then $\varphi$ is a group
monomorphism if and only if $B\varphi$ is a homotopy monomorphism.
\end{lemma}

The next lemma has an immediate algebraic proof, which is left to
the reader.
\begin{lemma} Let \mbox{$f\negmedspace : X \to Y$} and \mbox{$g\negmedspace : Y \to Z$}
be maps of spaces.
\begin{itemize}
  \item[(i)] If $f$ and $g$ are homotopy monomorphisms then the
composite \mbox{$g \circ f$} is a homotopy monomorphism.
  \item[(ii)] If $g \circ f$ is a homotopy monomorphism then $f$
  is a homotopy monomorphism.
\end{itemize}
\end{lemma}

We now arrive at the main purpose of this section.
\begin{definition}
A \emph{$p$-subgroup of a space $X$} is a pair $(P,f)$ where $P$
is a finite $p$-group and \mbox{$f \negmedspace : BP \to X$} is a
homotopy monomorphism. A $p$-subgroup $(S,f)$ of a space $X$ is a
\emph{Sylow $p$-subgroup of $X$} if every map $g\negmedspace : BP
\to X$, where $P$ is a finite $p$-group, factors up to homotopy
through \mbox{$f \negmedspace : BS \to X$}.
\end{definition}

\begin{example}
Let $G$ be a finite group with Sylow $p$-subgroup $S$ and let
$\iota_S$ be the inclusion \mbox{$S \hookrightarrow G$}. By
Sylow's second theorem, any homomorphism $\varphi \negmedspace : P
\to G$ factors through $\iota_S$ up to conjugacy. Recalling that
the classifying space functor induces a bijection
$$\Rep{P}{G} \stackrel{\cong}{\longrightarrow} [BP,BG],$$
we see that $(S,B\iota_S)$ is a Sylow $p$-subgroup of $BG$. Since
the $p$-completion functor induces a bijection
$$[BP,BG] \stackrel{\cong}{\longrightarrow} [BP,\pComp{BG}],$$
(this is a well known result, a proof of which can be found for
example in \cite{BL},) we see that $(S,\pComp{(B\iota_S)})$ is a
Sylow $p$-subgroup of $\pComp{BG}$.
\end{example}

The following lemma will be needed in the proof of the main
theorem of this paper and also to explain the difference between
the two stable classifications, as well as the difference between
the stable and unstable classifications.
\begin{lemma} \label{lem:SylowUnique}
Let $X$ be a space and let $(S,f)$ and $(S',f')$ be two Sylow
subgroups of $X$. Then there exists an isomorphism \mbox{$\gamma
\negmedspace : S \stackrel{\cong}{\longrightarrow} S'$} such that
\mbox{$f' \circ B\gamma \simeq f$}.
\end{lemma}
\begin{proof}
By the Sylow property of $(S',f')$, there is a map \mbox{$h
\negmedspace : BS \to BS'$} such that \mbox{$f' \circ h \simeq
f$}. Since $f$ is a homotopy monomorphism, $h$ must be a homotopy
monomorphism. Pick a \mbox{$\gamma \in \Hom{S}{S'}{}$} such that
\mbox{$h \simeq B\gamma $}. Then $\gamma$ is a group monomorphism
\mbox{$S \to S'$} such that \mbox{$f' \circ B\gamma \simeq f$}.
Similarly we get a monomorphism \mbox{$S' \to S.$} We deduce that
$S$ and $S'$ have the same order and therefore $\gamma$ is an
isomorphism.
\end{proof}

\section{Alternative stable classification} \label{sec:NewStableClassification}
In this section we prove the alternative stable classification
theorem stated in the introduction. We begin by recalling the
definition of fusion systems, which occur in the statement.
\begin{definition}
Let $G$ be a finite group and let $S$ be a Sylow subgroup of $G$.
The \emph{fusion system of $G$ over $S$} is the category $\F_S(G)$
whose objects are the subgroups of $S$, and whose morphisms are
the morphisms induced by conjugation in $G$:
$$\Hom{P}{Q}{\F_S(G)} = \Hom{P}{Q}{G}.$$

An isomorphism of fusion systems $\F_{S}(G)$ and $\F_{S'}(G)$ is a
group isomorphism \mbox{$\gamma \negmedspace : S
\stackrel{\cong}{\longrightarrow} S'$} such that for all
homomorphisms \mbox{$\varphi \negmedspace : P \to Q$} between
subgroups of $S$, we have \mbox{$\varphi \in \F_S(G)$} if and only
if \mbox{$\gamma \circ \varphi \circ \gamma^{-1} \in
\F_{S'}(G')$}.
\end{definition}
Martino-Priddy showed in \cite{MP} that the fusion system
$\F_S(G)$ can be recovered from the inclusion \mbox{$BS \to
\pComp{BG}$} via a simple homotopy theoretic construction. The
following proposition shows that the same construction works even
after infinite suspensions.
\begin{proposition}\label{prop:iotaToF}
Let $G$ be a finite group and let $S$ be a Sylow subgroup of $G$.
For subgroups $P$ and $Q$ of $S$ we have
$$\Hom{P}{Q}{\F_S(G)} = \{\varphi \in \Hom{P}{Q}{} \mid \ClSpectrum{\iota_Q} \circ \ClSpectrum{\varphi} \simeq  \ClSpectrum{\iota_P}\negmedspace : \ClSpectrum{P} \to \pComp{\ClSpectrum{G}}\}, $$
where $\iota_P$ and $\iota_Q$ denote the inclusions of  $P$ and
$Q$ in $G$.
\end{proposition}
\begin{proof}
Let \mbox{$\varphi \in \Hom{P}{Q}{}$}. By the Segal Conjecture,
\mbox{$ \ClSpectrum{\iota_Q} \circ \ClSpectrum{\varphi} \simeq
\ClSpectrum{\iota_P}$} if and only if \mbox{$[P,\iota_Q \circ
\varphi] = [P,\iota_P]$} in $\pComp{\widetilde{A}(P,G)}$, which is
clearly the case if and only if \mbox{$[P,\iota_Q \circ \varphi] =
[P,\iota_P]$} in $A(P,G)$. By definition, the last equality means
that there exist \mbox{$g \in P$} and \mbox{$h \in G$} such that
the following diagram commutes:
\[
\begin{CD}
{P} @ > {\iota_Q \circ \varphi} >> {G} \\
@ V \cong V c_g V @ VV c_h V \\
{P} @> {\iota_P} >> {G,} \\
\end{CD}
\]
or in other words such that
$$\varphi(x) = c_{h^{-1}g}(x)$$
for all \mbox{$x \in P$}. This is in turn true if and only
$\varphi$ is induced by a conjugation in $G$.
\end{proof}

Using Proposition \ref{prop:iotaToF}, and Oliver's solution of the
Martino-Priddy conjecture, the alternative stable classification
theorem stated in the introduction follows easily.
\begin{proof}[Proof of alternative stable classification]
Put
$$\F_1 := \F_{S}(G), \hspace{0.3 cm} \text{and} \hspace{.3cm} \F_2 := \F_{S'}(G').$$
$(i) \Rightarrow (ii):$ We show that $\gamma$ is an isomorphism of
fusion systems. For subgroups \mbox{$P, Q \leq S$}, and a
homomorphism \mbox{$\varphi \in \Hom{P}{Q}{\F_1}$}, we have
\begin{align*}
  \ClSpectrum{\iota_{\gamma(Q)}} \circ \ClSpectrum{(\gamma\vert_Q \circ \varphi \circ \gamma^{-1}\vert_{\gamma(P)})}
 &~{\simeq}~ h \circ \ClSpectrum{\iota_Q} \circ  \ClSpectrum{\varphi} \circ \ClSpectrum{\gamma^{-1}\vert_{\gamma(P)}}\\
 &~{\simeq}~ h \circ \ClSpectrum{\iota_P} \circ \ClSpectrum{\gamma^{-1}\vert_{\gamma(P)}}\\
 &~{\simeq}~ \ClSpectrum{\gamma\vert_P} \circ \ClSpectrum{\gamma^{-1}\vert_{\gamma(P)}}\\
 &~{\simeq}~ \ClSpectrum{\iota_{\gamma(P)}},
\end{align*}
so
$$\gamma\vert_Q \circ \varphi \circ \gamma^{-1}\vert_{\gamma(P)} \in \Hom{\gamma(P)}{\gamma(Q)}{\F_2}.$$
By symmetry we see that if
  \mbox{$\gamma\vert_Q \circ \varphi \circ \gamma^{-1}\vert_{\gamma(P)} \in \Hom{\gamma(P)}{\gamma(Q)}{\F_2},$}
then
  \mbox{$\varphi \in \Hom{P}{Q}{\F_1}$}. \\
$(ii) \Rightarrow (iii):$ This follows from
Oliver's proof of the
Martino-Priddy conjecture \cite{Ol1,Ol2}.\\
$(iii) \Rightarrow (i):$ Since $(S',B\iota_{S'})$ and $(S,h \circ
B\iota_{S})$ are both Sylow $p$-subgroups of $\pComp{BG'}$, there
is by Lemma \ref{lem:SylowUnique} an isomorphism \mbox{$\gamma
\negmedspace : S \to S'$} making the following diagram commute up
to homotopy:
  \[
\begin{CD}
{BS} @ > {B{\iota_{S}}} >> \pComp{B{G}} \\
@ V {B{\gamma}} VV @ VV h V \\
{B{S'}} @> {B{\iota_{S'}}} >> \pComp{B{G'}}. \\
\end{CD}
\]
$(i)$ follows upon infinite suspension.
\end{proof}

\section{Comparison of stable classifications} \label{sec:Discussion}
It is interesting to compare these two stable classifications. It
is easy to see that Condition $(i)$ in the alternative stable
classification implies Condition $(1)$ in the Martino-Priddy
classification stated in the introduction, and that Condition
$(ii)$ implies Conditions $(2)$ and $(3)$. However, the reverse
implications are not true. In \cite[Example 5.2]{MP3}
Martino-Priddy construct groups $G$ and $G'$, whose $p$-completed
classifying spaces are equivalent stably, but not unstably and it
is not difficult to see that their induced fusion systems are
non-isomorphic. Therefore the alternative classification theorem
is neither stronger nor weaker than the Martino-Priddy
classification theorem (in the sense that one cannot be deduced
from the other), but it does offer a `finer' classification of
$p$-completed stable classifying spaces of finite groups by
keeping track of more structure.

Recall that for a finite group $G$ with Sylow subgroup $S$, a
simple transfer argument shows that $\pComp{\ClSpectrum{G}}$ is a
wedge summand of $\ClSpectrum{S}$. In their classification,
Martino-Priddy regard $\pComp{\ClSpectrum{G}}$ as a wedge sum of
indecomposable stable summands of $\ClSpectrum{S}$. Thus two
finite groups $G$ and $G'$, both with Sylow subgroups isomorphic
to a finite $p$-group $S$, have stably homotopy equivalent
$p$-completed classifying spaces if and only if the multiplicity
of each indecomposable wedge summand of $\ClSpectrum{S}$ is the
same in $\pComp{\ClSpectrum{G}}$ as in $\pComp{\ClSpectrum{G'}}$.
Our new point of view is to differentiate between the summands of
$\ClSpectrum{S}$ and not just consider their homotopy types. By
adding the data of the maps \mbox{${\ClSpectrum{S}} \to
\pComp{\ClSpectrum{G}}$} and \mbox{$ \ClSpectrum{S} \to
\pComp{\ClSpectrum{G'}}$}, we obtain information about how the
summands of $\pComp{\ClSpectrum{G}}$ and $\pComp{\ClSpectrum{G'}}$
`sit inside' $\ClSpectrum{S}$. Instead of comparing the
multiplicities of homotopy types of summands, we compare whether
$\pComp{\ClSpectrum{G}}$ contains the same summands of
$\ClSpectrum{S}$ as $\pComp{\ClSpectrum{G'}}$ does.

The price of differentiating between summands is to restrict the
allowable maps between stable classifying spaces. For finite
groups $G$ and $G'$, with respective Sylow subgroups $S$ and $S'$,
any map \mbox{$h\negmedspace: \pComp{\ClSpectrum{G}} \to
\pComp{\ClSpectrum{G'}}$} can be extended to a map
\mbox{$\bar{h}\negmedspace: \ClSpectrum{S} \to \ClSpectrum{S'}$}
making the following diagram commute up to homotopy:
\[
\begin{CD}
{\ClSpectrum{S}} @ >>> {\pComp{\ClSpectrum{G}}} \\
@ V {\bar{h}} VV @ VV h V \\
{\ClSpectrum{S'}} @>>> {\pComp{\ClSpectrum{G'}}.} \\
\end{CD}
\]
The reason is that $\pComp{\ClSpectrum{G}}$ is a wedge summand of
$\ClSpectrum{S}$. It is by demanding in Condition $(i)$ that the
extension $\bar{h}$ be induced by a group homomorphism \mbox{$S
\to S'$} that we restrict the allowable homotopy types of $h$.
This way homotopy equivalent summands of $\ClSpectrum{S}$ can play
different roles, based on the group structure of $S$. It could be
interesting to reexamine the complete stable splittings of
classifying spaces of finite groups in \cite{BF} and \cite{MP2}
from this point of view, and attempt to label the stable summands
of the classifying space of a finite $p$-group in a way which
takes into account their role with respect to maps induced by
group homomorphisms.

%
%
It is also interesting that in Condition (iii) of the alternative
stable classification theorem, we apparently do not regard
$\pComp{BG}$ and $\pComp{BG'}$ as spaces under their Sylow
subgroups. By Lemma \ref{lem:SylowUnique} the homotopy type of
$\pComp{BG}$ determines not only the isomorphism class of its
Sylow subgroup $S$, but also the homotopy type of the pair
$(BS,B\iota_S)$ as an object over $\pComp{BG}$. In particular the
homotopy type of the arrow \mbox{$(BS \xrightarrow{B\iota}
\pComp{BG})$} is implicit in the homotopy type of $\pComp{BG}$.

In the stable setting, Nishida proved in \cite{Ni} that the stable
homotopy type of $\pComp{BG}$ determines the isomorphism class of
the Sylow subgroup $S$. Since $\pComp{\ClSpectrum{G}}$ is a wedge
summand of $\ClSpectrum{S}$, and the stable splitting of $BS$ is
unique up to order and homotopy equivalence of summands
(\cite{MP2,Ni,He}), there is a homotopy unique spectrum $X$ such
that \mbox{$\ClSpectrum{S} \simeq \pComp{\ClSpectrum{G}} \vee X$},
and we have
$$(\ClSpectrum{S}\xrightarrow{\ClSpectrum{\iota}} \pComp{\ClSpectrum{G}}) \simeq (\pComp{\ClSpectrum{G}} \vee X \xrightarrow{proj} \pComp{\ClSpectrum{G}}).$$
Therefore the homotopy type of the arrow
\mbox{$(\ClSpectrum{S}\xrightarrow{\ClSpectrum{\iota}}
\pComp{\ClSpectrum{G}})$} is also implicit in the homotopy type of
$\pComp{\ClSpectrum{G}}$.

The difference between the stable and unstable cases is that in
the unstable case an equivalence between classifying spaces of
Sylow subgroups $(BS,B\iota_S)$ and $(BS',B\iota_{S'})$ is a
homotopy equivalence \mbox{$h \negmedspace : BS
\xrightarrow{\simeq} BS'$} respecting the inclusion maps (up to
homotopy) and this map must be induced by a group isomorphism. In
the stable case, an equivalence between stable classifying spaces
of Sylow subgroups $(\ClSpectrum{S},\ClSpectrum{\iota_S})$ and
$(\ClSpectrum{S},\ClSpectrum{\iota_S})$ is similarly a homotopy
equivalence \mbox{$h \negmedspace : \ClSpectrum{S}
\xrightarrow{\simeq} \ClSpectrum{S'}$} respecting the inclusion
maps (up to homotopy), but this map need not at all be induced by
a group isomorphism (although by \cite{Ni} its existence does
imply that the groups are isomorphic). Therefore we need to impose
a condition on the allowable maps between stable classifying space
of finite $p$-groups which is automatically satisfied for the
unstable classifying spaces.


\begin{thebibliography}{9}

\bibitem{BF} D. Benson, M. Feshbach, \emph{Stable splittings of
classifying spaces of finite groups}, Topology \textbf{31} (1992),
157-176.

\bibitem{Monster} A.K. Bousfield, D.M. Kan, \emph{Homotopy limits, completions and
localizations}, Lecture Notes in Mathematics, vol. \textbf{304},
Springer Verlag, 1972.

\bibitem{BL} C. Broto, R. Levi, \emph{On spaces of self homotopy equivalences of $p$-completed classifying spaces of finite groups and homotopy group
extensions}, Topology \textbf{41} (2002), 229-255.


\bibitem{Car} G. Carlsson, \emph{Equivariant stable homotopy and Segal's Burnside ring conjecture,}
Ann. of Math. \textbf{120} (1984), 189-224.


\bibitem{Ev} L. Evens, \emph{The cohomology ring of a finite
group,} Trans. Amer. Math. Soc. \textbf{101} (1961), 224-239.

\bibitem{He} H.-W. Henn, \emph{Some finiteness results in the category of unstable modules over the Steenrod algebra and applications to stable
splittings,} Forschungschwerpunlet Geometrie Univ. Heidelberg
\textbf{47} (1989).

\bibitem{LMM} L.G. Lewis, J.P. May, J.E. McClure, \emph{Classifying $G$-spaces and the Segal conjecture,}
Current trends in algebraic topology, Part 2 (London, Ont., 1981)
CMS Conf. Proc. 2, (1981), 165-179.

\bibitem{MP2}
 J. Martino, S. Priddy, \emph{The complete stable splitting for the classifying space of a finite
 group}, Topology \textbf{31} (1992), 143-156.

\bibitem{MP3}
J. Martino, S. Priddy, \emph{Stable homotopy classification of
 $\pComp{BG}$}, Topology \textbf{34} (1995), 633-649.

\bibitem{MP}
 J. Martino, S. Priddy, \emph{Unstable homotopy classification of
 $\pComp{BG}$}, Math. Proc. Cambridge Phil. Soc. \textbf{119}
 (1996), 119-137.

\bibitem{Ni} G. Nishida \emph{Stable homotopy type of classifying
spaces of finite groups}, Algebraic and Topological Theories
(1985), 391-404.

\bibitem{Ol1}
 B. Oliver, \emph{Equivalences of classifying spaces completed at odd primes,}
 Math. Proc. Camb. Phil. Soc. \textbf{137} (2004), 321-347.

\bibitem{Ol2}
 B. Oliver, \emph{Equivalences of classifying spaces completed at the prime two},
 Mem. Amer. Math. Soc. (To appear).

\end{thebibliography}
\end{document}